\documentclass[12pt,draft]{article}
\usepackage{latexsym,color,amsmath,amsthm,amssymb,amscd,amsfonts}

\newtheorem{lemma}{Lemma}[section]

\newtheorem{theorem}[lemma]{Theorem}
\newtheorem{definition}[lemma]{Definition}
\newtheorem{corollary}[lemma]{Corollary}

\newtheorem{example}[lemma]{Example}

\newtheorem*{remark*}{Remark}

\def\C{C}

\begin{document}
\title {The Nonlinear
Maslov index and the Calabi homomorphism}
\author{Gabi Ben Simon
\footnote{This paper is a part of the author's Ph.D. thesis, being
carried out under the supervision of Professor Leonid Polterovich,
at Tel-Aviv University.} \footnote{Partially supported by Israel
Science Foundation grant 11/03.}
\\
School of Mathematical Sciences \\
Tel Aviv University
69978 Tel Aviv, Israel \\
 gabiben@post.tau.ac.il} \maketitle
\begin{abstract} \noindent In this paper we find that the
 asymptotic nonlinear Maslov index defined on the universal cover of the group of all contact Hamiltonian diffeomorphisms
of the standard $(2n-1)$-dimensional contact sphere is a
quasimorphism. Then we show our main result: Let $M$ be standard
$(n-1)$-dimensional complex projective space. We prove that the value of the
pullback of the asymptotic nonlinear Maslov index to the
universal cover of the group of Hamiltonian diffeomorphisms of
$M$, when evaluated on a diffeomorphism supported in a sufficiently small open
subset of $M$, equals  $\frac{-n}{2\pi vol(M)}$ times the Calabi
invariant of this diffeomorphism.
%restricted to elements supported in a
% sufficiently small open subset of the complex projective space. We
%get that it
  %, where $vol(M)$ stands for the standard volume of the
%complex projective space.

\end{abstract}

\section{Introduction and results}

\subsection{Basics on the nonlinear Maslov index} \label{basics on maslov}

 Let $(S^{2n-1},\xi)$ be
 the $(2n-1)$-dimensional sphere endowed with the standard contact
 structure. We denote by $\mathcal{C}$ the
 identity component of the group of those contact Hamiltonian diffeomorphisms of $(S^{2n-1},\xi)$, which commute with the inversion $x \mapsto -x$ on $S^{2n-1}$
  (see~\cite{M-S} for preliminaries on contact manifolds).
 Let $\mathcal{\widetilde{C}}$ be the
 universal cover of $\mathcal{C}$. Recall that an element $\tilde{f}$ of $\mathcal{\widetilde{C}}$
can be identified as the time-1 map of the flow generated  by a
contact Hamiltonian, up to homotopy with fixed endpoints.
Multiplication on $\mathcal{\widetilde{C}}$ is defined as follows:
let $\tilde{f}$, $\tilde{g}$ in $\mathcal{\widetilde{C}}$ be the
time-1 map of the flows $\{\tilde{f}_{t}\}_{0\leq t \leq 1}$ and
$\{\tilde{g}_{t}\}_{0\leq t \leq 1}$ respectively. Then the
multiplication of $\tilde{f}$ and $\tilde{g}$ is the homotopy
class of the path $\{\tilde{f}_{t}\tilde{g}_{t}\}_{0\leq t \leq
1}.$

Moreover, $\mathcal{C}$ carries the structure of an infinite
dimensional manifold. Following Givental (see ~\cite{G}), we
denote by $\Delta(\mathcal{C})$ the \textit{discriminant} of
$\mathcal{C}$, that is the hypersurface of $\mathcal{C}$ formed by
all contactomorphisms which have a fixed point, and preserve the
standard contact form at that point. Using generating functions,
Givental has defined a co-orientation of $\Delta(\mathcal{C})$ in
$\mathcal{C}$.

Given a path $\{\tilde{f}_{t}\}_{t\geq 0}$ in $\mathcal{C}$, an
intersection index can be assigned to it by counting all signed
intersections of the path with a subspace of
$\Delta(\mathcal{C})$, called the $\it{train}$ and denoted by
$\Gamma(\mathcal{C})$, relative to the co-orientation at
nonsingular points of $\Delta(\mathcal{C})$, (see again~\cite{G}).

We call the intersection index of $\{\tilde{f}_{s}\}_{s\in[0,t]}$
with $\Gamma(\mathcal{C})$  \textit{the nonlinear Maslov index}
of $\{\tilde{f}_{s}\}_{s\in[0,t]}$ and denote it by $\mu(t)$.
Following Givental we now define the \textit{asymptotic nonlinear
Maslov index}: let $\tilde{f}$ be an element of $\mathcal{\widetilde{C}}$ generated
by a time-1 periodic contact Hamiltonian $\widetilde{F}:
S^{2n-1}\times S^{1} \longrightarrow \mathbb{R}$. Let
$\{\tilde{f}_{t}\}_{t\geq 0}$ be the flow generated by
$\widetilde{F}$, which is a path in $\mathcal{C}$. The asymptotic nonlinear Maslov index of $\tilde{f}$ is then given by
\begin{equation}\label{maslov for universal}\textsl{m}(\tilde{f}):=\lim\limits_{t\rightarrow\infty}
\frac{\mu(t)}{t}\,.\end{equation}
%$\begin{equation} \textsl{m}(\tilde{f})=\textsl{m}(\tilde{f}_{t})
%\lim\limits_{t\rightarrow\infty}
%\frac{\mu(\tilde{f}_{t})}{t}
%$\end{equation}.

\subsection{The quasimorphism property}\label{quasi property}
\begin{definition} \label{quasi}
{\rm Let $G$ be a group. A \textit{quasimorphism} $r$ on $G$ is a
function $r:G\rightarrow\mathbb{R}$ which satisfies the
homomorphism equation up to a bounded error: there exists $R>0$
such that
$$|r(fg)-r(f)-r(g)|\leq R$$
for all $f,g\in G$}.
\end{definition}
See ~\cite{Ba} for preliminaries on quasimorphisms. A
quasimorphism $r_{h}$ is called \textit{homogeneous} if
$r_{h}(g^{m})=mr_{h}(g)$ for all $g\in G$ and $m\in \mathbb{Z}$.
Every quasimorphism $r$ gives rise to a homogeneous one
$$r_{h}(g)=\lim \limits_{m\rightarrow\infty}\frac{r(g^{m})}{m}.$$

 We now can state the first theorem of this paper
which is actually an easy consequence of ~\cite{G}.
\begin{theorem} \label{m is quasi} The asymptotic nonlinear Maslov index
$\textsl{m}:\mathcal{\widetilde{C}}\longrightarrow\mathbb{R}$ is a
quasimorphism.
\end{theorem} The proof will be given in Section 2 after some
preliminaries on the asymptotic nonlinear Maslov index.

\subsection{The pullback of \textsl{m}}\label{pull back of m}Let
$(\mathbb{C}P^{n-1},\omega)$ be $(n-1)$-dimensional complex
projective space endowed with the Fubini-Study form, normalized so that $\int\limits_{\mathbb{C}P^{1}}\omega=\pi.$ Let $G=Ham(\mathbb{C}P^{n-1},\omega)$ be the identity
component of the group of all Hamiltonian diffeomorphisms of
$(\mathbb{C}P^{n-1},\omega)$  and $\widetilde{G}$ its universal
cover (see~\cite{M-S} and~\cite{P} for preliminaries on $G$). Let
$\{f_{t}\}_{t\in[0,1]}$ be a path in $G$ representing an element
of $\widetilde{G}$. There is a unique Hamiltonian $F:\mathbb{C}P^{n-1}\times
S^{1}\longrightarrow\mathbb{R}$, which generates $\{f_{t}\}_{t\in[0,1]}$ and satisfies the normalization
condition:
$$\int\limits_{\mathbb{C}P^{n-1}}{F_{t}\omega^{n-1}}=0$$
for all $t\in S^{1}$.

\noindent We denote the space of all normalized $1$-periodic Hamiltonians by
$\textsf{F}$. In order to pull back the asymptotic nonlinear
Maslov index from $\mathcal{\widetilde{C}}$ to  $\widetilde{G}$ we
will use the Hopf bundle $$ p:S^{2n-1}\longrightarrow
\mathbb{C}P^{n-1}\text{   .}$$
\noindent Let $f \in \widetilde{G}$, and let $\{f_{t}\}_{t\in[0,1]}$ be a representative of $f$. Let $F\in\textsf{F}$ be the unique normalized Hamiltonian generating this path. Then we can lift $F$ to
$S^{2n-1}$ to a time dependent function,
$$\widetilde{F}_{t}:=p^{\ast}F_{t}\,.$$
%satisfies the condition:
%$$\widetilde{F}(\tilde{x}_{0})=F(x_{0})\text{   }\forall\tilde{x}_{0}\in \pi^{-1}_{1}(x_{0})$$
%For every $x_{0}\in\mathbb{C}P^{n-1}$.
$\widetilde{F}$ is constant along the fibers of the fiber bundle
$p :S^{2n-1}\longrightarrow \mathbb{C}P^{n-1}$. This $\widetilde{F}$ is a contact Hamiltonian on the contact sphere. The contact flow
$\{\tilde{f}_{t}\}_{t\in[0,1]}$ generated by this contact
Hamiltonian constitutes a lift of the flow
$\{f_{t}\}_{t\in[0,1]}$, and represents an element $\tilde f$ of
$\mathcal{\widetilde{C}}$. The value of the pullback of the nonlinear Maslov index at \(f\) is then set to be $m(\tilde f)$.
We will denote by \begin{equation} i:
\widetilde{G}\longrightarrow \mathcal{\widetilde{C}}\end{equation}
the above injection, and by
\begin{equation}i^{\ast}\textsl{m}:
\widetilde{G}\longrightarrow\mathbb{R}\end{equation}
the pullback of the asymptotic nonlinear Maslov index.
As a corollary of Theorem~\ref{m is quasi} we have

\begin{corollary} The pullback of the asymptotic nonlinear Maslov index $$i^{\ast}\textsl{m}:
\widetilde{G}\longrightarrow\mathbb{R}$$ is a quasimorphism.
\end{corollary}

\textbf{Remark.} Any Hamiltonian (not necessarily normalized)
generating an element of $\widetilde{G}$ can be lifted via $p$ in
the way described above . The normalization condition is required
since we want to embed $\widetilde{G}$ in
$\mathcal{\widetilde{C}}$.
%Hamiltonian which is a lift of a
%Hamiltonian $F$, constitute a subgroup of $\mathcal{C}$ called the
%\textit{quantomorphisms} group.

We close this subsection with a simple, but important example.
%Note that the lifted Hamiltonians are .
\begin{example}{\rm Consider the flow $f_{t}=\mathbf{1}$ for all
$t$ where $ \mathbf{1}$ is the identity map of
$\mathbb{C}P^{n-1}$. Any function of the form $F(x,t)=F(t)$
(clearly, not necessarily normalized)
 generates this flow. Let $p^{\ast}F_{t}$ be its lift to the
 sphere. An easy computation shows that the flow $\{\tilde{f}_{t}\}$,
 generated by $p^{\ast}F_{t}$,
  leaves the fibers of the Hopf bundle invariant, and  rotates
them by angle $e^{iF(t)}$ for time $t$.}
\end{example}

\subsection{The Calabi property of $i^{\ast}$\textsl{m}}
\label{calabi of i*m}

   In this subsection we will assume that $G=Ham(M,\omega)$, where
$M^{2n}$ is any closed connected symplectic manifold.

\noindent The group $G$ has a natural class of subgroups $G_{U}$
associated to nonempty open subsets $U\subset M$. An element $f$ of $G$
lies in $G_{U}$ if it is the time-$1$ map of a Hamiltonian
flow $\{f_{t}\}_{t\in[0,1]}$ with $f_{0}=\mathbf{1}$, generated by a
Hamiltonian $F_{t}:M\longrightarrow\mathbb{R}, \; t\in [0,1]$, such that
\begin{equation}\text{support}\,(F_{t})\subset U \text{ for all }t.  \label{support condition}\end{equation}

 Consider the map $Cal_{U}:G_{U}\longrightarrow\mathbb{R}$ given
 by\begin{equation} \label{cal integral}
 f \mapsto \int_{0}^{1}dt \int_{M}
 F_{t}\omega^{n}\,.
 \end{equation}
 When the symplectic form $\omega$ is exact on $U$, this map is
 well-defined, meaning that it does not depend on the specific
 choice of the Hamiltonian $F$ generating $f$. In fact, $Cal_{U}$
 is a homomorphism called \textit{the Calabi homomorphism}
 (see~\cite{B},~\cite{C} and~\cite{W} for further information).

\noindent The above definitions can be ``lifted'' to
$\widetilde{G}$: \noindent for a nonempty subset $U\subset M$,
consider the subgroup
 $\widetilde{G}_U\subset \widetilde{G}$ defined as follows. An
 element of  $\widetilde{G}$ lies in $\widetilde{G}_U$ if it can be represented by a Hamiltonian flow
 $\{f_{t}\}_{t\in[0,1]}$ with $f_{0}=\mathbf{1}$, generated by a
 Hamiltonian $F_{t}$ satisfying condition \eqref{support
 condition}. Formula \eqref{cal integral} gives rise to a well-defined homomorphism
 $$\widetilde{Cal}_U:\widetilde{G}_U\longrightarrow\mathbb{R}\,.$$
 In what follows we deal with the class $\mathcal{D}$ of all nonempty open subsets $U$ which can be \textit{displaced} by a
 Hamiltonian diffeomorphism: \begin{equation} hU\cap
 \text{Closure} (U)=\varnothing  \text{ for some }  h\in G\text{
 .}
 \end{equation}
 %Put $$\mathcal{D}_{ex}=\{U\in \mathcal{D}| \omega\text{ is exact
 %on }U\} \text{ . }$$

  We need one more definition before stating our main
  result.

\begin{definition}{\rm A quasimorphism on $\widetilde{G}$ coinciding
with the Calabi homomorphism
$\{\widetilde{Cal}_U:\widetilde{G}_U \rightarrow
 \mathbb{R}\}$ on each $U \in \mathcal{D}$ is said to have
the \textit{Calabi property}. } \end{definition} From now on
$M=\mathbb{C}P^{n-1}$. Here is our main result:

\begin{theorem} \label{main theorem} The pullback of the asymptotic nonlinear Maslov index
$i^{\ast}\textsl{m}: \widetilde{G}\longrightarrow\mathbb{R}$  has
the Calabi property up to multiplication by $-\frac{n}{2\pi
\textup{vol}\,(M)}$. That is
$$i^{\ast}\textsl{m}([f]) = -\frac{n}{2\pi\textup{vol}\,(M)}\, \widetilde{Cal}_U([f])$$
for every element $[f]\in\widetilde{G}_{U}$.
\end{theorem}
The proof will be given in section 3.

$\bf{Remark.}$ Given a Hamiltonian $H:M \times
S^{1}\longrightarrow \mathbb{R}$ denote by
$f^H$ the  element of $\widetilde{G}$
generated by $H$. A function
  $r:\widetilde{G}\longrightarrow\mathbb{R}$ is called \textit{continuous}, if it satisfies the
  following condition: if a sequence $\{H_{i}\}$ of smooth
Hamiltonians $H_{i} : M \times S^{1} \longrightarrow \mathbb{R}$ $C^{0}$-converges
 to a smooth function $H:M \times S^{1}\longrightarrow \mathbb{R}$ then
  \begin{equation} r(f^{H_{i}})\longrightarrow
  r(f^{H})\text{ as }i \rightarrow \infty.\end{equation}

\noindent It follows from Corollary 3 of Section 9 of~\cite{G}
that the pullback of the asymptotic nonlinear Maslov index is
continuous. Summarizing, the pullback of the asymptotic
nonlinear Maslov index $i^{\ast}\textsl{m}:
\widetilde{G}\longrightarrow\mathbb{R}$ is a continuous
quasimorphism with the Calabi property (up to multiplication by a constant).

\subsection{Other quasimorphisms with the Calabi
property}Another construction of a  continuous quasimorphism on
$\widetilde{G}$ with the Calabi property has been carried out by
Entov and Polterovich in~\cite{En-P}. In particular it applies to the manifold $(\mathbb{C}P^{n-1},\omega)$.
%is spherically monotone, with a semi
%simple quantum homology algebra and thus has the above
%quasimorphism with the Calabi property.

Their construction uses the spectral invariants of elements of
$\widetilde{G}$. These spectral invariants are defined via the Floer theory on the symplectic manifold. The first
step in our proof of theorem~\ref{main theorem} goes along the lines
of~\cite{En-P}. Moreover, step 3 of the proof of~\ref{main
theorem} was partly motivated by the inspection of the behavior
of spectral invariants defined in terms of generating functions
(see~\cite{En-P} for more details).

\bigskip

\section{Preliminaries on the nonlinear Maslov index}
In this short section we state the main properties of the
asymptotic nonlinear Maslov index and the nonlinear Maslov index
$\mu$, which we use in the proof of the main result
(theorem~\ref{main theorem}).

First, we remind the reader that the sign of an intersection of a
path in $\widetilde{\mathcal C}$ with the train (which is a subspace of the discriminant) is
determined via the co-orientation of the train, which is defined by methods of
generating functions. The two basic facts that we need are
grouped together in the following theorem, which is due to Givental
(see~\cite{G}, ~\cite{G1}, and~\cite{T}).

\begin{theorem}~\label{G theorem} (a) The
nonlinear Maslov index (relative to the train) of any two paths in
$\mathcal{C}$ is the same, provided that there is a homotopy
between them, such that the end points of the paths in
this homotopy lie outside the discriminant.

 (b) Let $f$ be a Hamiltonian map
 of $\mathbb{C}P^{n-1}$ and let $\tilde{f}$ be its
 lift to $\mathcal{C}$. Assume that it has no fixed points.
 Then the nonlinear Maslov index
 of the loop $\{e^{\pm2\pi
 it}\tilde{f}\}_{t\in[0,1]}$ (relative to the train) is $\pm n$.

\end{theorem}

\subsection{Proof of Theorem ~\ref{m is quasi}}
Let $\tilde{f}\in\widetilde{\mathcal C}$. By definition, it is generated by
a $1$-periodic contact Hamiltonian. Let $\{\tilde{f}_{t}\}$ be
the contact flow generated by this Hamiltonian. Then we have
$\tilde{f}^{m}=\tilde{f}_{m}$ for all integer $m$. This means that
\begin{equation}~\label{with m}\textsl{m}(\tilde{f})=
\lim\limits_{m\rightarrow\infty}\frac{\mu(\tilde{f}^{m})}{m}\,.\end{equation}

It suffices to prove that
$\mu$ is a quasimorphism on $\widetilde{\mathcal{C}}$ (see section~\ref{quasi property}).
For this, let $[\tilde{f}]$ and $[\tilde{h}]$ be elements of $\widetilde{\mathcal{C}}$, represented by paths $\{\tilde{f}_{t}\}_{t\in[0,1]}$ and $\{\tilde{h}_{t}\}_{t\in[0,1]}$, respectively. Consider the following
path $\gamma$, which represents $[\tilde{f}][\tilde{h}]$. It is the
concatenation
$$\{\tilde{f}_{2t}\}_{0\leq t \leq \frac{1}{2}} \ast
\{\tilde{f}_{1}\tilde{h}_{2t-1}\}_{\frac{1}{2}\leq t\leq
1}:=\begin{cases}\tilde{f}_{2t} &0 \leq t \leq \frac{1}{2}\\
\tilde{f}_{1}\tilde{h}_{2t-1} & \frac{1}{2} \leq t \leq 1
\end{cases} \;.$$

It follows from Givental's construction that $\mu$ is additive
with respect to  concatenation. So we have that
\begin{equation}\label{catenation}  \mu(\gamma)= \mu([\{\tilde{f}_{2t}\}])+
\mu([\{\tilde{f}_{1}\tilde{h}_{2t-1}\}]).\end{equation}
Moreover it is also follows from his construction that
\begin{equation}\label{mu inequality} |\mu([\{\tilde{h}_{2t-1}\}])
-\mu([\{\tilde{f}_{1}\tilde{h}_{2t-1}\}])|\leq 2n . \end{equation}
Combining \eqref{catenation} and \eqref{mu inequality} we get
$$|\mu(\gamma)-\mu([\{\tilde{f}_{2t}\}])-\mu([\{\tilde{h}_{2t-1}\}])|\leq 2n .$$
Since the paths $\{\tilde{f}_{2t}\}$ and  $\{\tilde{h}_{2t-1}\}$
represent respectively $[\tilde{f}]$ and $[\tilde{h}]$, we have:
$$|\mu([\tilde{f}][\tilde{h}])-\mu([\tilde{f}])-\mu([\tilde{h}])|\leq 2n \,, $$
which is the desired. \qed

\section{Proof of Theorem ~\ref{main theorem}}
\subsection{Step 1}

Let $\{f_{t}\}_{t\in[0,1]}$ be a representative of an element
$[f] \in \widetilde{G}_{U}$ for some open set $U$ in $\mathcal{D}$,
generated by a time-1 periodic Hamiltonian $F_{t}$  (see
~\ref{calabi of i*m}). We have to show that:
\begin{equation}\label{the equation} -\frac{n}{2\pi \textup{vol}\,(M)} \, \widetilde{Cal}_U([f]) =
i^{*}\textsl{m}([f])\;.
\end{equation}

Recall that the right hand side is actually the asymptotic nonlinear Maslov
index of the lift $\{\tilde{f}_{t}\}_{t\geq 0}$ of the flow $\{f_{t}\}_{t\geq0}$
(see~\ref{basics on maslov} and~\ref{pull back of m}). Let $h$ be a Hamiltonian diffeomorphism displacing $U$,
 and suppose that it is the time-$1$ map of a Hamiltonian isotopy $\{h_t\}_{t \in [0,1]}$. As a first
step, we claim that it is enough to show that the asymptotic
nonlinear Maslov index of the path $\{\widetilde{hf_{t}}\}_{t\geq
0}$, which is the lift of the path $\{hf_{t}\}_{t\geq 0}$, equals the left-hand side of
~\eqref{the equation}. This follows easily from inequality
~\eqref{mu inequality}. Since the difference of the
nonlinear Maslov indexes of these paths is bounded by $2n$, they have the same asymptotic nonlinear Maslov index.
%To show this, it is
%enough to show that
%$\textsl{m}(\widetilde{hf}_{t})=\textsl{m}(\tilde{f}_{t})$. By
%definition $h$ is the time-1 map of a Hamiltonian flow
%$\{h_{t}\}_{t\in[0,1]}$. Let $\{\tilde{h}_{t}\}_{t\in [0,1]}$ be
%\its lift, and  consider the concatenation
%$\{\tilde{h}_{t}\}_{t\in[0,1]}\ast \{\widetilde{hf}_{t-1}\}_{t\geq
%1}$. By the additivity of $\mu$ with respect to concatenation we
%get:
%$$\textsl{m}(\{\widetilde{hf_{t}}\})=\textsl{m}(\{\tilde{f}_{t}\})$$
%As required.

From now on we will focus on the concatenation
$$\{{h}_{t}\}_{t\in[0,1]}\ast\{{hf_{t-1}}\}_{ t\geq1}:=\gamma_{t}\,,$$
and its lift to $\mathcal{C}$
$$\{\tilde{h}_{t}\}_{t\in[0,1]}\ast \{\widetilde{hf}_{t-1}\}_{t\geq1}:
=\widetilde{\gamma}_{t}\,.$$

%\textbf{CHECK HERE}

Let us emphasize, that the paths $\{\gamma_{t}\}_{0\leq t\leq
t_{0}}$  and $\{\widetilde{\gamma}_{t}\}_{0\leq t \leq t_{0}}$
(for $1 \leq t_{0}$) represent respectively $hf_{t_{0}-1}$ and
$\widetilde{hf}_{t_{0}-1}$ as elements of $\widetilde{G}$ and
$\widetilde{\mathcal{C}}$ respectively.

Secondly, note that because $h$ displaces $U$, the set of fixed points of $\{hf_{t}\}_{t\geq0}$
coincides with the set of fixed points of $h$ for all $t\geq0$.
Let $x_{0}$ be a fixed point of
$h$, and let $\widetilde{x_{0}}\in p^{-1}(x_{0})$ be a point of
the fiber $p^{-1}(x_{0})$ above $x_{0}$. Of course
$\widetilde{\gamma}_{t}(\widetilde{x}_0)$ is a lift of
$\gamma_{t}(x_{0})$. In  particular for $0\leq t \leq 1$ the path
$\widetilde{\gamma}_{t}(\widetilde{x}_0)$ is a lift of the loop
$\gamma_{t}(x_{0})$, but this path need not be closed. Nevertheless
the point $\widetilde{x}_{0}$ will return to its fiber for $t=1$.
For $t\geq1$ we have $\gamma_{t}(x_{0})=x_{0}$, so the fiber
$p^{-1}(x_{0})$ is invariant under the flow
$\widetilde{\gamma}_{t}$. Therefore we have:
\begin{equation}\label{rotation equation}\widetilde{\gamma}_{t}(\widetilde{x_{0}})=
e^{i\theta(t)}\widetilde{x_{0}}\text{  ,  }t\geq 1\text{  }
\end{equation}
where $\theta(t)$ is determined by the contact flow  generated by
the lift of the Hamiltonian generating the path
$\gamma_{t}$.
\subsection{Step 2}

Now, assume that for $t=s$, $s>1$ the path
$\widetilde{\gamma}_{t}$ intersects the train $\Gamma(\mathcal{C})$. This
means that there exists $x_{0}\in\{\text{fix}\,(h)\}$ such that in
equation \eqref{rotation equation} we have $\theta(s)=2 \pi k$ for
some $k\in\mathbb{Z}$. In the second step of the proof we first
calculate the function $\theta(s)$, and then modify the path
$\widetilde{\gamma}_{t}$ to one for which the function
$\theta(s)$ is linear. This will help us to make
the final reduction step in the next section, where we
actually calculate the asymptotic nonlinear Maslov index.
%is to show  that in a small
%neighborhood of $s$ $\text{sign}(\widetilde{Cal}([f]))$
%determines the direction of rotation of the fiber
%$p^{-1}(x_{0})$ (this will ,lead us to the fact that
%$\text{sign}(\widetilde{Cal}([f]))$  determines the sign of the co-orientation of
%$\widetilde{\gamma}_{t}$ each time it intersect the train).

\noindent \textbf{Notation:} From now on we denote
$C = -\frac{\widetilde{Cal}_U([f])}{\textup{vol}\,(M)}$.

Let $\{K_{t}\}$ be the normalized Hamiltonian generating the flow $\{\gamma_{t}\}$, and let  $\{\widetilde{K}_{t}\}$
be its lift. Clearly for $t\geq 1$ we have:
$$K_{t}=F_{t}+c(t)$$
%\frac{1}{vol{M}} \int \limits_{M}F_{t}\omega^{n-1}$$
where $F_{t}$ is a $1$-periodic Hamiltonian generating the flow $\{f_{t}\}$, and
$$c(t)=-\frac{1}{\textup{vol}\,{M}} \int\limits_{M}F_{t}\omega^{n-1}\;.$$
Moreover, for all $x\in M\setminus U$ we have $K_{t}(x)=c(t)$.

For $s\geq 1$ define
$$\theta(s):=\int_{1}^{s}c(t)\,dt=-\int_{1}^{s}dt\:\frac{1}{\textup{vol}\,{M}} \int \limits_{M}F_{t}\omega^{n-1}\;.$$
Since $F_{t}$ is periodic we have for every $m\in \mathbb{N}$ and $0\leq v \leq 1$
$$\theta(m+v)=C(m-1) - \int_{0}^{v}dt\:\frac{1}{\textup{vol}\,{M}} \int
\limits_{M}F_{t}\omega^{n-1}\;.$$

Note also that the function $\varphi(s):= \theta(s)-C(s-1)$ is $1$-periodic
and satisfies
$\varphi(m)=0$ for every $1\leq m\in\mathbb{N}$. In particular for all such $m$ we have
\begin{equation}\label{m}C(m-1)=\theta(m)\,.\end{equation}

The notation $\theta(s)$ here is deliberate. Actually, it is the solution of \eqref{rotation equation} for all fixed points of $h$.
% So according to inequality~\ref{L} and equality~\ref{m} the average number of signed
%nonlinear Maslovs of $\{\widetilde{\gamma}_{t}\}$ with the train will not change if we substitute
%$$\theta(s)=\C(s-1)+\theta_{j}$$ where $0\leq\theta_{j}<2\pi$
%is defined by the condition $$e^{i \theta_{j}}(\tilde{x}_{j})=\tilde{x}_{j}$$
% for all $\tilde{x}_{j}\in p^{-1}(x_{j})$ at time $s=1$.The average number of times in which
% $\{\widetilde{\gamma}_{t}\}$ hit the train is the asymptotic nonlinear maslov index.
Note that $\theta(s)$ measures the total displacement of the fiber after time $s$. From now on we assume, without loss of generality, that $h$ has a finite number of fixed points. This always can be achieved by a small perturbation of $h$.
%We say that $x_{0}\in\{\text{fix}(h)\}$ \textit{occurs} at $t=s_{0}$ if $\widetilde{\gamma}_{s_{0}}$ is
%the identity map when restricted to the fiber $p^{-1}(x_{0})$.

Now let $\{x_{1}, \dots, x_{l} \}$ be the set of fixed points of
$h$ and let $\tilde{x}_{j}$ be a point above $x_{j}$ for all
$1\leq j \leq l $. Let $\alpha_{j}$, $0\leq \alpha_{j}<2\pi$, be
the angle defined by the equation
$\widetilde{\gamma}_{1}(\tilde{x}_{j})=e^{i\alpha_{j}}\tilde{x}_{j}$.
We have for $s\geq1$
\begin{equation} ~\label{rotation with phi} \widetilde{\gamma}_{s} (\tilde{x}_{j})=e^{i\C(s-1)}e^{i\varphi(s)}e^{i\alpha_{j}}(\tilde{x}_{j})\,. \end{equation}

    From the properties of $\varphi(s)$ we get that the path
    $\{e^{i\varphi(s_{0}+t)}\}_{t\in[0,k]}$, $k\in \mathbb{N}$
     is homotopic to the constant loop $\{e^{i\varphi(s_{0})}\}$, for every $s_{0}\geq 0 $.
     We emphasize that the homotopy has $\{e^{i\varphi(s_{0})}\}$
   as the starting and end point of all of its paths.

   Let us denote such a homotopy by $\{e^{iF_{s}(t)}\}_{(s,t) \in
   [0,1]\times [0,k]}$. So we have $$F_{0}(t)=
   \varphi(s_{0}+t),\text{  }F_{1}(t)= \varphi(s_{0})$$ and
$$F_{s}(0)=F_{s}(k)=\varphi(s_{0}).$$
   For reasons that will become clear we are
   looking for $s_{0}$ that satisfies the following condition:
   \begin{equation}\C(s_{0}+k-1)+\varphi(s_{0}+k)+\alpha_{j}
   \notin 2\pi \mathbb{Z}\label{choose_s_0} \end{equation}
for all $1\leq j \leq l$ and $k\in \mathbb{N}$.

   Such $s_{0}$ clearly exists. To see this, first note that the set $2\pi \mathbb{Z}$ is
   countable as well as $\mathbb{N}$, the set $\{\alpha_{1},\dots,\alpha_{l}\}$ is finite,
   and finally $\varphi(s_{0}+k)=\varphi(s_{0})$
   for all $k \in \mathbb{N}$. So we have to
   find a point in the image of the function $\C(s-1)+
   \varphi(s)$ which belongs to the complement of the countable set
   $\bigcup \limits_{j,k}\{2\pi \mathbb{Z}- \alpha_{j} - k\C\}$. This function is continuous and unbounded,
   thus its image clearly contains a segment, and consequently a point in the complement of the above set.

   For such $s_{0}$ we have for all $j,k$    \begin{equation} ~\label{rotation with phi,k}
    \widetilde{\gamma}_{s_{0}+k}(\tilde{x}_{j})= e^{i\C(s_{0}+k-1)}e^{i\varphi(s_{0})}
    e^{i\alpha_{j}}(\tilde{x}_{j})\neq \tilde{x}_j\,. \end{equation}

    That is for all $k \in \mathbb{N}$\; $\widetilde{\gamma}_{s_{0}+k}$ has no fixed
    points (and thus it is outside of the train). Now we have come to the key idea of step 2.

    \begin{lemma} The family of maps
    $$\psi(s,t)=\{\widetilde{\gamma}_{s_{0}+t}e^{-iF_{s}(t)}e^{i
    \varphi(s_{0})}\}_{(s,t) \in [0,1]\times [0,k]}$$ is a homotopy of paths in $\mathcal{C}$ which has the
    following properties:

    (a) For $s=1$ it is the path
    $\{\widetilde{\gamma}_{s_{0}+t}\}_{t \in [0,k]}$ .

(b) The starting point and the end point of each one of the paths
comprising this homotopy has no fixed points. That is, the maps
$\psi(s,0)$, $\psi(s,k)$ have no fixed points.\end{lemma}
\begin{proof} For $s=1$ we have $F_{1}(t)= \varphi(s_{0})$. Thus by the very definition of $\psi$ we get part (a).
For $t=0$ and $t=k$ we know that for every $0\leq s \leq k$
$F_{s}(0)=F_{s}(k)=\varphi(s_{0})$. Thus  by the definition
of $\psi$ and formula \eqref{rotation with phi,k}  we get part (b)
of the lemma.
\end{proof}

It follows from the invariance of the nonlinear Maslov index with
respect to homotopies (see theorem~\ref{G theorem}) that the path
$\psi(0,t)$ has the same index as the path
 $\psi(1,t)$ (which is
$\{\widetilde{\gamma}_{s_{0}+t}\}_{t \in [0,k]}$). Now,
restricting the path $\psi(0,t)$  to the fibers over the fixed
points we have the following equation
\begin{equation} ~\label{rotation with residue}
 \psi(0,t)(\tilde{x}_{j})=
    e^{i\C(s_{0}+t-1)}e^{i\varphi(s_{0})}e^{i
    \alpha_{j}}\tilde{x}_{j}\,.\end{equation}

    We get the following equation by restricting the path
    $e^{i[C-\varphi(s_{0})]}\psi(0,t)$ to the fibers over
    the fixed points:$$e^{i[C-\varphi(s_{0})]}\psi(0,t)(\tilde{x}_{j})=e^{i\C(s_{0}+t)}e^{i
    \alpha_{j}}(\tilde{x}_{j})$$

Moreover, we see from equation~\eqref{mu inequality} that the
difference between the nonlinear Maslov indexes of
%Recall that we have shown in the proof of theorem~\ref{m is
    %quasi} that $\mu$ has the quasimorphism property. From this we
    %get that the asymptotic nonlinear Maslov index of the paths
    $e^{i[C-\varphi(s_{0})]}\psi(0,t)$ and $\psi(0,t)$ is bounded by $2n$ for all $k\in
    \mathbb{N}$. Thus, by the same argument we have used in step 1, we
    obtain that the asymptotic nonlinear Maslov indexes of the paths $e^{i[C-\varphi(s_{0})]}\psi(0,t)$ and
    $\psi(0,t)$, $t\geq0$, are equal.

   %$\mathbf{(I)}$  For some $m\in \mathbb{N} $ we have  $\widetilde{\gamma}_{m}(\tilde{x}_{j})\neq\tilde{x}_{j}$
    %for all $1\leq j \leq k $

  %  $\mathbf{ (II) }$ $\C\notin 2\pi \mathbb{Q} $

  It follows from the above discussion that without loss of
  generality we can assume that for $s\geq s_{0}$, $s_{0}\geq1$
\begin{equation}\label{cal is rotation}\widetilde{\gamma}_{s}(\tilde{x}_{j})=
e^{i[\C(s-s_{0})+\theta_{j}]}\tilde{x}_{j}\,,\end{equation}
where $\theta_{j}$ is the \noindent initial condition, which depends on
$s_{0}$ and $x_{j}$, satisfies $0\leq \theta_{j}< 2\pi$ and is defined by the
equation
$$\widetilde{\gamma}_{s_{0}}(\tilde{x}_{j})=e^{i\theta_{j}}\tilde{x}_{j}\,.$$

\subsection{Step 3}

%We start with the following lemma:

%\begin{lemma}\label{occurs in integer value} A point $x_{0}\in\{fix(h)\}$
%occurs at $t=s_{0}$, $s_{0}>1$ of the path
%$\{\widetilde{\gamma}_{t}\}$ if and only if the action of $x_{0}$
%(as a fixed point of $\{hf_{s_{0}}\}$) is $0$ modulo
%$2\pi$.\end{lemma}

%\noindent The meaning of this lemma is that each point
%$x_{0}\in\{fix(h)\}$ is "detected" by the flow
%$\{\widetilde{\gamma}_{t}\}$ each time its action has an integer mod
%$2\pi$ value.

%The proof of lemma~\ref{occurs in integer value} is essentially the
%same as the proof of lemma~\ref{action as integral of alpha} only
%that now we use equality~\ref{v} (see~\cite{T} lemma 5.8 for a
%detailed proof).

Before introducing the last step of the proof let us motivate it.
We have seen that we may assume that the fibers above the fixed
points of $h$ rotate, under the action of $\{\widetilde\gamma_t\}$, at the same angular speed and in the same
direction of rotation (where the orientation is induced by the
action of the group $\{e^{it}\}$). So after a period of time of length
$\frac{2\pi}{|C|}$ each of the above fibers makes full rotation.
This might lead us to believe that the number of signed intersections of $\{\widetilde\gamma_t\}$, where $t$ runs from some $t_0$ to $t_0 + \frac{2\pi}{|C|}$, with the train, is independent of $t_0$, provided that $\widetilde\gamma_{t_0}$ lies outside the train. In order to show that this is indeed the situation, we need the following lemma.

\begin{lemma}\label{second homotopy} Let
$\{\widetilde{\gamma}_{t}\}_{ t\geq1}$ be the path in
$\mathcal{C}$ as defined at the beginning of this section, and assume that it satisfies \eqref{cal is rotation}. Choose
$s_{0}$ such that $\widetilde{\gamma}_{s_{0}}$ is without fixed
points. Then the family of maps
$$\{e^{\C i(1-s)t}
\widetilde{\gamma}_{s_{0}+st}\}
_{(s,t)\in[0,1]\times[0,\frac{2\pi}{|\C|}]}$$ is a homotopy
between the paths
  $\{\widetilde{\gamma}_{s_{0}+t}\}_{t\in[0,\frac{2\pi}{|\C|}]}$    and
  $\{e^{\C it}\widetilde{\gamma}_{s_{0}}\}_{t\in[0,\frac{2\pi}{|\C|}]}$
  such that the starting and the end points of the paths of this homotopy stay outside the discriminant.
\end{lemma}

\begin{proof}The first part of the lemma is trivial, since all we
need to do is to substitute $s=0$  and  $s=1$ to get the above
paths. It is also obvious that the homotopy stays in $\mathcal{C}$.
What we must show is that the constant path
$$\Lambda_{1}^{s}:=\{\widetilde{\gamma}_{s_{0}}\}$$
and the path
$$\Lambda_{2}^{s}:=\{e^{\C i(1-s)\frac{2\pi}{|\C|}}
\widetilde{\gamma}_{s_{0}+\frac{2\pi s}{|\C|}}\}_{s\in[0,1]}$$
have no fixed points.
We need only check this for the path $\Lambda_{2}$. To this end let $x_{0} \in \text{fix}(h)$, and let $\tilde{x}_{0}$ be any point
above $x_{0}$. Let $0<\theta_{0}<2\pi$ be the angle of rotation of
$\tilde{x}_{0}$ at time $s_{0}$, where $s_0$ is chosen so as to satisfy \eqref{choose_s_0}.
According to \eqref{cal is rotation} we have
$$\widetilde{\gamma}_{s_{0}+\frac{2\pi s}{|\C|}}(\tilde{x}_{0})
=e^{i(sign(\C)2\pi s+\theta_{0})}(\tilde{x}_{0})\,.$$
From this we get
 $$\Lambda_{2}^{s}(\tilde{x}_{0})=
 e^{i[sign(\C)(1-s)2\pi+sign(\C)2\pi s+\theta_{0}]}(\tilde{x}_{0})\,,$$
 or simply:
 $$\Lambda_{2}^{s}(\tilde{x}_{0})=e^{i(sign(\C)2\pi+\theta_{0})}(\tilde{x}_{0})\,.$$
 According to our choice of $\theta_{0}$ clearly
 $\Lambda_{2}^{s}(\tilde{x}_{0})\neq\tilde{x}_{0}$ for all
 $s\in[0,1]$. This concludes the proof of our lemma.
 \end{proof}

 %Now, It is proved in
 %~\cite{G} that any two paths homotopic to each
 %other, such that the ends of the paths of this homotopy  stays
 %outside of the train, has the same nonlinear Maslov index with the
 % train.

Here we appeal to part $(a)$ of theorem~\ref{G theorem} and to lemma~\ref{second
homotopy} in order to conclude that the paths  $$\{\widetilde{\gamma}_{s_{0}+t}\}_{t\in[0,\frac{2\pi}{|\C|}]}$$    and
  $$\{e^{\C it}\widetilde{\gamma}_{s_{0}}\}_{t\in[0,\frac{2\pi}{|\C|}]}$$
  has the same nonlinear Maslov index relative to the train. The calculation of  the nonlinear Maslov
  index of $\{e^{\C it}\widetilde{\gamma}_{s_{0}}\}_{t\in[0,\frac{2\pi}{|\C|}]}$
  %can ,easily, be obtained using~\cite{T}. This is due to the
  %following two facts proved in~\cite{T} :
uses part $(b)$ of theorem~\ref{G theorem}, from which we
get
\begin{corollary}\label{loop index} The nonlinear Maslov
index of the path
$\{\widetilde{\gamma}_{s_{0}+t}\}_{t\in[0,\frac{2\pi}{|\C|}]}$
relative to the train is $\textup{sign}(C)n$.\end{corollary}

% \noindent \textbf{Fact 2.} If there exist $t_{0}\in[0,1)$ such that
 %the co-homological length of $M^{-}(t)$, of the path  $\{e^{2\pi
 %it}\widetilde{\zeta}_{s_{0}}\}_{t\in[0,1]}$, changes by more then $1$
 %while passing through $t_{0}$, then $\zeta_{s_{0}}$ has infinite
 %number of fixed points.

 %Now, by our construction, the path $\{\gamma_{t}\}$ has finite
 %number of fixed points for all $t$'s .

% co-homological length of $M^{-}(t)$

  %$\{e^{\Cit}\widetilde{\gamma}_{s_{0}}\}_{t\in[0,\frac{2\pi}{|\C|}]}$
  %changes each time it intersect the train by 1.
   %This with fact 1 gives

 %that the path
  %$\{e^{\Cit}\widetilde{\gamma}_{s_{0}}\}_{t\in[0,\frac{2\pi}{|\C|}]}$
%intersect the train exactly $sign(\C)n$ times. We conclude that
%the path

Now we can calculate the asymptotic nonlinear Maslov index of $[f]$. Choose
$s_{0}\geq 1$ such that $\widetilde{\gamma}_{s_{0}}$ has no fixed points. Note that
$\widetilde{\gamma}_{s_{0}+\frac{2\pi k}{|\C|}}$ also has
no fixed points, and thus is outside the train, for all $k\in
\mathbb{N}$. Define
\begin{align*}a_{k}&:=\mu\Big(\frac{2\pi k}{|C|}+s_{0}\Big)-\mu(s_{0})\\
&= \sum\limits_{d=1}^{k}\Big[\mu\Big(\frac{2\pi
d}{|C|}+s_{0}\Big)-\mu\Big(\frac{2\pi
(d-1)}{|C|}+s_{0}\Big)\Big]\,.
\end{align*}
According to corollary~\ref{loop index} this equals
$$kn\, \text{sign}\,(C).$$ Thus
\begin{align*}i^{\ast}\textsl{m}([f])&=\lim\limits_{t\rightarrow\infty}\frac{\mu(t)}{t}=\lim\limits_{k\rightarrow\infty}\frac{a_{k}}{\frac{2\pi
k}{|C|}}\\
&=\frac{n\,\text{sign}\,(C)|C|}{2\pi}=\frac{nC}{2\pi}\\
&=-\frac{n}{2\pi \textup{vol}\,(M)}\,\widetilde{Cal}_U([f])\,,
\end{align*}
  %|\C|})}{\frac{2\pi k}{|\C|}}=
  %\lim\limits_{k\rightarrow\infty}\frac{\mu(\frac{2\pi
%k}{|\C|}+s_{0})-\mu(s_{0})} {\frac{2\pi k} {|\C|}}=$$
%$$\lim\limits_{k\rightarrow\infty}\frac{\mu(\frac{2\pi
%k}{|\C|}+s_{0})-\mu(\frac{2\pi (k-1)}{|\C|}+s_{0}) +\mu(\frac{2\pi
%(k-1)}{|\C|}+s_{0})-...-\mu(s_{0})}{\frac{2\pi k} {|\C|}}=$$
%$$\lim\limits_{k\rightarrow\infty}\frac{|\C|}{2\pi}[\frac{\mu(\frac{2\pi
%k}{|\C|}+s_{0})-\mu(\frac{2\pi
%(k-1)}{|\C|}+s_{0})}{k}+...+\frac{\mu(\frac{2\pi
% }{|\C|}+s_{0})-\mu(s_{0})}{k}]$$
 %$$=\lim\limits_{k\rightarrow\infty}\frac{|\C|}{2\pi}[\frac{sign(\C)nk}{k}]=\frac{n\C}{2\pi}$$
as claimed.

Finally, the case $\C=0$ follows from the considerations of step
2. According to what we have seen, we can take a homotopy of
the path $\widetilde{\gamma}_{s_{0}+t}$ such that there occurs no rotation in
the fibers (see \eqref{rotation with residue}), so clearly the asymptotic nonlinear Maslov index of $[f]$ is zero.
% according to this formula rotation of the fibers remains
%bounded. And clearly the asymptotic nonlinear Maslov index is zero.

\bigskip

\noindent$\bf{Acknowledgments.}$ My sincere gratitude goes to
Leonid Polterovich for all his support, ideas, and encouragement.
I want to thank Michael Entov for reading drafts of this paper,
and for coming to my talk about this work at Tel Aviv university
and sharing his opinions with me. I would also like to thank Yaron
Ostrover and Frol Zapolsky  for reading the last draft of this
paper and for their remarks.

\bigskip

\noindent

\end{document}